\documentclass[a4,12pt]{article}

\usepackage{amsmath, amssymb, amsthm, multicol}

\usepackage[english]{babel}

\newtheorem{lemma}{Lemma}
\newtheorem{theorem}{Theorem}

\newtheorem{prop}{Proposition}
\newtheorem{coro}{Corollary}

\allowdisplaybreaks

\def\N{{\mathbb N}} 
\def\Z{{\mathbb Z}}

\def\R{{\mathbb R}} 
\def\C{{\mathbb C}}

\def\s{{\bf s}}

\def\x{{\bf x}}
\def\y{{\bf y}}
\def\z{{\bf z}}
\def\a{{\bf a}}
\def\b{{\bf b}}
\def\k{{\bf k}}

\def\X{{\bf X}}

\def\u{{\bf u}}

\def\q{{\bf q}}

\def\Nb{{\bf N}}
\def\P{{\bf P}}

\newcommand{\zerob}{\boldsymbol{0}}

\newcommand{\alphab}{\boldsymbol{\alpha}}
\newcommand{\betab}{\boldsymbol{\beta}}

\newcommand{\gammab}{\boldsymbol{\gamma}}
\newcommand{\thetab}{\boldsymbol{\theta}}
\newcommand{\omegab}{\boldsymbol{\omega}}

\newcommand{\mub}{\boldsymbol{\mu}}
\newcommand{\deltab}{\boldsymbol{\delta}}

\def\eps{{ \varepsilon }}
\newcommand{\be}{\begin{enumerate}}
\newcommand{\ee}{\end{enumerate}}
\let\ds=\displaystyle

\usepackage{xcolor}

\begin{document}
%%***********************************************************
\title{Values at non-positive integers of generalized Euler-Zagier multiple zeta-functions\footnote{The authors benefit from the financial support of 
the French-Japanese Project ``Zeta-functions of Several Variables and  Applications" 
(PRC CNRS/JSPS 2015-2016).}}
\author{Driss Essouabri \hskip 0.4cm and \hskip 0.4cm
Kohji Matsumoto}
\date{}
\maketitle 
%\vspace{3mm}

\noindent
{{\bf {Abstract.}}\\
{\it {\small We give new closed and explicit formulas for ``Multiple zeta values'' at non-positive integers of generalized Euler-Zagier multiple zeta-functions.
%and its partially twisted analogues. In the non-twisted case, 
We first prove these formulas for a small convenient class of these multiple zeta-functions and then use the analyticity of the values on the parameters defining the multiple zeta-functions to deduce the formulas in the general case. Also, for our aim we prove an extension of     "Raabe's lemma" due to E. Friedman and A. Pereira (Lemma 2.4 of \cite{fredman}). 
%In the partially twisted case, we use a different  method based on M. de Crisenoy's result \cite{dC} on the fully twisted case and the Mellin-Barnes integral formula.
}}

\medskip

\noindent
{\small {\bf Mathematics Subject Classifications: Primary 11M32; Secondary 11M41.} }
%\vskip .1 in
\setcounter{tocdepth}{2}
% \tableofcontents
%\vspace{15mm}

%%%%%%%%%%%%%%%%%%%%%%%%%%%%%%%%%%%%%%%%%%%%%%%%%%%%%%%%%%%%%%%%%%%%%%%%%%%%%%%%%%%%%%%%%
\section {Introduction and the statement of main results}
%%%%%%%%%%%%%%%%%%%%%%%%%%%%%%%%%%%%%%%%%%%%%%%%%%%%%%%%%%%%%%%%%%%%%%%%%%%%%%%%%%%%%%%%%

Let $\mathbb{N}$, $\mathbb{N}_0$, $\mathbb{Z}$, $\mathbb{R}$, and $\mathbb{C}$ be the sets of 
positive integers, non-negative integers, 
rational integers, real numbers, and complex numbers, respectively.

Let $\gammab =(\gamma_1,\dots, \gamma_n) \in \C^n$ and  $\b =(b_1,\dots, b_n) \in \C^n$ be  such that 
$\Re (\gamma_j) >0$ and $\Re (b_j) >-\Re (\gamma_1)$ for all $j=1,\dots, n$.

The generalized  Euler-Zagier multiple zeta-function is defined formally 
for  $n-$tuples of complex variables $\s=(s_1,\dots, s_n)$ by
\begin{equation}\label{EZMzetadef}
\zeta_n(\s; \gammab; \b) :=
\sum_{m_1\geq 1 \atop m_2, \dots, m_n \geq 0} \frac{1}{\prod_{j=1}^n (\gamma_1 m_1+\dots+ \gamma_j m_j +b_j)^{s_j}}.
\end{equation}
If $b_1=0$, $b_j=j-1$ for all $j=2,\dots, n$ 
and $\gamma_j=1$ for all $j=1,\dots, n$, then $\zeta_n(\s; \gammab; \b)$ coincides with 
the classical  Euler-Zagier multiple zeta-function (see \cite{zagier} and \cite{hoffman})
$$
\sum_{1\leq m_1< m_2<\dots <m_n } \frac{1}{m_1^{s_1}\dots m_n^{s_n}}.
$$
The generalized Euler-Zagier multiple zeta-function $\zeta_n(\s; \gammab; \b)$ converges absolutely in the domain 
\begin{equation}\label{domaincv}
\mathcal D_n:=\{\s=(s_1,\dots, s_n)\in \C^n \mid \Re (s_j+\dots+s_n) >n+1-j~~{\rm for \; all}\; j=1,\dots, n\}
\end{equation}
(see \cite{MatIllinois}),
and has a meromorphic continuation to $\C^n$ whose poles are located in the union of the hyperplanes 
 $$s_j+\dots+s_n =(n+1-j)-k_j \quad (1\leq j\leq n,~k_1,\dots, k_n \in \N_0).$$
Moreover, it is known that for $n\geq 2$, almost all  $n-$tuples of non-positive integers lie on the singular locus above and are points of discontinuity (see \cite{AET}, Th.1).
The evaluation of (limit) values of multiple zeta-functions at those points was first considered by
S. Akiyama, S. Egami and Y. Tanigawa \cite{AET}, and then studied further 
by \cite{AT}, \cite{sasaki1},
\cite{sasaki2}, \cite{komori}, \cite{onozuka}, and \cite{MOW}.

In \cite{komori}, Y. Komori proved that for  $\Nb =(N_1,\dots, N_n)\in \N_0^n$ and $\thetab=(\theta_1,\dots, \theta_n) \in \C^n$ such that $\theta_j+\dots+\theta_n \neq 0$ for all $j=1,\dots, n$, 
the limit 
\begin{equation}\label{gmzvtheta}
\zeta_n^{\thetab} (-\Nb; \gammab; \b):=\lim_{t\rightarrow 0}\zeta_n(-\Nb+t\thetab; \gammab; \b)
\end{equation} 
exists, and expressed it in terms of $\Nb$, $\thetab$ and generalized multiple Bernoulli numbers defined implicitly as coefficients of some multiple series. \par
Our main result (i.e. Theorem \ref{main1}) gives a closed explicit formula for $\zeta_n^{\thetab} (-\Nb; \gammab; \b)$ in terms of $\Nb$, $\thetab$ and {\it only classical Bernoulli numbers}
$B_k$ $(k\in \N_0)$ defined by
\begin{equation}\label{bernoulli}
\frac{x}{e^x-1}=\sum_{k=0}^{\infty}B_k\frac{x^k}{k!}.
\end{equation}

Before giving our result let us introduce a few notations:
\be
\item For  $\mathbf{x}=(x_1,\ldots,x_n)\in \C^n$, we write $|\mathbf{x}|=x_1+\cdots+x_n$;
\item For $\mathbf{x}=(x_1,\ldots,x_n)\in \C^n$ and $\mathbf{k}=(k_1,\ldots,k_n)\in \N_0^n$, we write
$$\ds \mathbf{x}^{\mathbf{k}}=\prod_{i=1}^n x_i^{k_i},\qquad
{\mathbf{x} \choose \mathbf{k}}=\prod_{i=1}^n {x_i \choose k_i}{\rm ;}
$$
\item For $\Nb =(N_1,\dots, N_n) \in \N_0^n$ and $\alphab=(\alpha_1,\ldots,\alpha_n) \in \N_0^n$, we define 
\begin{equation}\label{kalphadef}
K(\Nb, \alphab):=\left\{j\in\{1,\dots,n\} \left| (n+1-j)+ \sum_{i=j}^n N_i
=\sum_{i=j}^n \alpha_i\right.\right\}
\end{equation}
and
\begin{equation}\label{ialphadef} 
L(\Nb, \alphab):=\big \{j\in \{1,\dots, n\}\mid \alpha_j \geq N_j+1\big\};
\end{equation}
\item For $\Nb =(N_1,\dots, N_n) \in \N_0^n$ and $I\subset \{1,\dots, n\}$, we define 
\begin{eqnarray}\label{jindef}
\mathcal J(I,\Nb):=\{\alphab \in \N_0^n \mid K(\Nb,\alphab)=I
 {\mbox { and }} |L(\Nb, \alphab)|= |I|\};
\end{eqnarray}
{\bf Remark:} $\mathcal J(I, \Nb)$ is a finite set and $\mathcal J(I,N) \subset \{0,\dots, |\Nb|+n\}^n$.
(See Lemma \ref{finitesetlemma} for a proof of this fact).
\item 
For $\alphab=(\alpha_1,\ldots,\alpha_n) \in \N_0^n$ and  $\b=(b_1,\ldots,b_n) \in \C^n$ 
we define the polynomial (in $\b$)  $c_n(\b;\alphab, \k)$ 
(where $\k=(k_1,\ldots,k_n) \in \N_0^n, |\k|\leq |\alphab|)$) as the coefficients of the polynomial 
$\prod_{j=1}^n (\sum_{i=1}^j X_i +b_j)^{\alpha_j}$; that is 
\begin{equation}\label{benoulliexpension}
\prod_{j=1}^n (\sum_{i=1}^j X_i +b_j)^{\alpha_j}= \sum_{\k \in \N_0^n,\atop  |\k|\leq |\alphab|} 
c_n(\b;\alphab, \k) ~\X^\k = \sum_{\k \in \N_0^n,\atop  |\k|\leq |\alphab|}
 c_n(\b;\alphab, \k) ~X_1^{k_1} \dots X_n^{k_n}.
 \end{equation}
\ee
With these notations our main result is the following:
\begin{theorem}\label{main1}
Let $\gammab =(\gamma_1,\dots, \gamma_n) \in \C^n$ and  $\b =(b_1,\dots, b_n) \in \C^n$ be  such that 
$\Re (\gamma_j) >0$ and $\Re (b_j) > -\Re (\gamma_1)$ for all $j=1,\dots, n$.
Let $\thetab=(\theta_1,\dots, \theta_n) \in \C^n$ be such that $\theta_j+\dots+\theta_n \neq 0$ for all $j=1,\dots, n$.
Then, for any $\Nb =(N_1,\dots, N_n) \in \N_0^n$,\\
 the limit 
$$\zeta_n^{\thetab} (-\Nb; \gammab; \b):=\lim_{t\rightarrow 0}\zeta_n(-\Nb+t\thetab; \gammab; \b)$$
exists, and is explicitly given by
\begin{eqnarray}\label{thm-1-formula}
&&\zeta_n^{\thetab} (-\Nb; \gammab; \b)\nonumber\\
&=&\sum_{I\subset \{1,\dots, n\}} \sum_{\alphab \in \mathcal J (I, \Nb)} \sum_{ \k \in \N_0^n,\atop  |\k|\leq |\alphab|}
  \frac{c_n(\b;\alphab, \k)~ (-1)^{n-|I|+\sum_{j\in L(\Nb, \alphab)} (\alpha_j -N_j)} \prod_{j \not \in L(\Nb, \alphab)} {N_j \choose \alpha_j}}{\prod_{j\in L(\Nb, \alphab)}
 \alpha_j {\alpha_j -1\choose N_j}
 \prod_{j\not \in I} 
 \left(n+1-j+\sum_{i=j}^n N_i-\sum_{i=j}^n\alpha_i\right)}\nonumber\\
 &&\qquad \times \left(\gamma_1^{|\Nb|-|\alphab|+n+k_1-1}\prod_{j=2}^n \gamma_j^{k_j-1}\right)~\left(\frac{ \prod_{j  \in I(\alphab, N)} \theta_j}{\prod_{j \in I}(\theta_j+\dots+\theta_n)}\right)~\left( \prod_{j=1}^n B_{k_j}\right).
\end{eqnarray}
\end{theorem}

An essential idea in our proof of Theorem \ref{main1} is to prove these formulas first for a small convenient class of these multiple zeta-functions and then use the analyticity of the values on the parameters defining the multiple zeta-functions to deduce the formulas in the general case.    We also prove an extension of     a lemma of "Raabe type" due to E. Friedman and A. Pereira (Lemma 2.4 of \cite{fredman}) and use it
in the proof. 

\section{Some useful lemmas}
%%%%%%%%%%%%%%%%%%%%%%%%%%%%%%%%%%%%%%%%%%%%%%%%%%%%%%%%%%%%%%%%%%%%%%%%%%%%%%%%%%%%%%%%%%
\begin{lemma}\label{integralexpression}
Let $\gammab =(\gamma_1,\dots, \gamma_n) \in \C^n$  be such that 
$\Re (\gamma_j) >0$ for any $j=1,\dots, n$.
Define for  $\s=(s_1,\dots, s_n) \in \mathcal D_n$ (see (\ref{domaincv}))
\begin{equation}\label{ynsgammaadef}
Y_n(\s; \gammab):=\int_{(1,\infty)\times (0,\infty)^{n-1}} \prod_{j=1}^n \left(\sum_{i=1}^j \gamma_i x_i\right)^{-s_j}~dx_n\dots dx_1.
\end{equation}
Then,  for  $\s=(s_1,\dots, s_n) \in \mathcal D_n$, $Y_n(\s; \gammab)$ is absolutely convergent and 
$$Y_n(\s; \gammab)=\frac{\gamma_1^{-s_1-\dots-s_n+n}}{(\gamma_1\dots \gamma_n) \prod_{j=1}^n (s_j+\dots+s_n +j-n-1)}.$$
In particular, $Y_n(\s; \gammab)$ has a meromorphic continuation to $\C^n$ and its polar locus is the set 
$$\bigcup_{j=1}^n \{\s=(s_1,\dots, s_n) \in \C^n\mid s_j+\dots +s_n =n+1-j \}.$$
\end{lemma}
{\bf Proof of Lemma \ref{integralexpression}:}   Just integrate first with respect to the variable $x_n$ and then with respect to $x_{n-1}$ etc.\qed

\begin{lemma}\label{finitesetlemma}
Let $\Nb =(N_1,\dots, N_n) \in \N_0^n$ and $I\subset \{1,\dots, n\}$. 
The set $\mathcal J(I,\Nb)$ defined by \eqref{jindef}
is a finite set and $\mathcal J(I,\Nb) \subset \{0,\dots, |\Nb|+n\}^n$.
\end{lemma}
{\bf Proof of Lemma \ref{finitesetlemma}:}\\
Denote by  $j_1,\dots, j_q$ the elements of the set $I$, where $q=|I|$.  We assume without loss of generality that $j_1<j_2<\dots<j_q$.

Let $\alphab=(\alpha_1,\ldots,\alpha_n) \in \mathcal J(I,\Nb)$. 
It follows that for any $k=2,\dots,q$,
\begin{eqnarray*}
\sum_{j=j_{k-1}}^{j_k-1} \alpha_j  &= & \sum_{j=j_{k-1}}^{n} \alpha_j - \sum_{j=j_{k}}^{n} \alpha_j\\
&=&(n+1-j_{k-1}) + \big(\sum_{j=j_{k-1}}^{n} N_j\big) -(n+1-j_k)- \big(\sum_{j=j_{k}}^{n} N_j\big) \\
&=& (j_k-j_{k-1})+\sum_{j=j_{k-1}}^{j_k-1} N_j \\
&\geq & 1+ \sum_{j=j_{k-1}}^{j_k-1} N_j,
\end{eqnarray*}
hence $[j_{k-1},  j_k) \cap L(\Nb, \alphab) \neq \varnothing$ for all $k=2,\dots,q$.
Moreover, the identity 
$$ \sum_{j=j_{q}}^{n} \alpha_j = (n+1-j_q)+\sum_{j=j_{q}}^{n} N_j  \geq  1+\sum_{j=j_{q}}^{n} N_j $$ 
implies also that $[j_q,  n] \cap L(\Nb, \alphab) \neq \varnothing$.
Since $| L(\Nb, \alphab)|= q$, the above observation implies that 
$ \min  L(\Nb, \alphab) \geq j_1$.
We deduce that for  $j\in  L(\Nb, \alphab)$,
$$\alpha_j \leq \sum_{j=j_{1}}^{n} \alpha_j = (n+1-j_1)+\sum_{j=j_{1}}^{n} N_j  \leq |\Nb|+n.$$
If $j\notin L(\Nb,\alphab)$, obviously
$\alpha_j<N_j+1\leq |\Nb|+n$.
This ends the proof of Lemma \ref{finitesetlemma}.\qed\par
The following lemma is crucial for our proof of Theorem \ref{main1}.
\begin{lemma}\label{keycomputation}
Let $\Nb=(N_1,\dots, N_n) \in \N_0^n$, $\alphab=(\alpha_1,\dots, \alpha_n)\in \N_0^n$ and $\thetab=(\theta_1,\dots, \theta_n) \in \C^n$ be such that 
$\theta_j+\dots+\theta_n \neq 0$ for all $j=1,\dots, n$.
Set 
$$\delta :=\frac{1}{2} \min\left\{(1+|\theta_j|)^{-1}, |\theta_j+\dots+\theta_n|^{-1};~j=1,\dots, n\right\} \in (0, 1/2).
$$
Let $U_\delta :=\{t\in \C;~ |t|<\delta\}$. Define for 
$t\in U_{\delta}\setminus\{0\}$:
\begin{equation}\label{def-G}
G_{\Nb, \alphab, \thetab}(t):= \frac{\prod_{j=1}^n {N_j-t\theta_j \choose \alpha_j} }{\prod_{j=1}^n\left(t-\frac{(N_j+\dots+N_n)+(n+1-j)-(\alpha_j+\dots+\alpha_n)}{\theta_j+\dots+\theta_n}\right)}.
\end{equation}
Let $q=q(\Nb, \alphab):=|K(\Nb, \alphab)|$ and $q'=q'(\Nb, \alphab):=|L(\Nb, \alphab)|$,
where $K(\Nb,\alphab)$, $L(\Nb,\alphab)$ are defined by \eqref{kalphadef}, \eqref{ialphadef},
respectively.
Then,
\be 
\item $q'\geq q$;
\item $G_{\Nb, \alphab, \thetab}(t)$ is {\it analytic} in the disk $U_{\delta}$
and there exists a constant $C=C(\Nb, \theta) >0$ (which is independent of $\alphab$)
such that 
$$|G_{\Nb, \alphab, \thetab}(t)| \leq C~|t|^{q'-q} \quad {\rm for\;all}\;\; t\in U_\delta.$$
\item If $q'> q$, then $G_{\Nb, \alphab, \thetab}(0)=0$;
\item If $q'=q$, then 
\be
\item
$$G_{\Nb, \alphab, \thetab}(0)=\frac{(-1)^{n-q}\left(\prod_{j\in L(\Nb, \alphab)} \frac{(-1)^{\alpha_j -N_j} \theta_j}{\alpha_j {\alpha_j -1 \choose N_j}}\right) 
\prod_{j\not \in  L(\Nb, \alphab)} {N_j \choose \alpha_j}}{\prod_{j\not \in K(\Nb, \alphab)}\left(\frac{(N_j+\dots+N_n)+(n+1-j)-(\alpha_j+\dots+\alpha_n)}{\theta_j+\dots+\theta_n}\right)};$$
\item $\alphab \in \mathcal J \left (K(\Nb, \alphab), N\right) \subset \{0,\dots, |\Nb|+n\}^n$ 
(see Lemma \ref{finitesetlemma}).
\ee
\ee 
\end{lemma}

{\bf Proof of Lemma \ref{keycomputation}:}\\
$\bullet$ {\bf Proof of point 1:} 
Repeat the argument of the proof of the previous lemma with $I=K(\Nb,\alphab)$.
It follows that $q'=| L(\Nb, \alphab)|\geq q$.\qed\\
$\bullet$ {\bf Proof of point 2:} 
First it is easy to see that $G_{\Nb, \alphab, \thetab}(t)$ is {\it analytic} in all the pointed disk  $U_\delta \setminus\{0\}$.
(Since $(N_j+\dots+N_n)+(n+1-j)-(\alpha_j+\dots+\alpha_n)\in\mathbb{Z}$, if it is not zero, then
$\left|\frac{(N_j+\dots+N_n)+(n+1-j)-(\alpha_j+\dots+\alpha_n)}{\theta_j+\cdots+\theta_n}\right|\geq 2\delta$.)
Moreover we have for $j=1,\dots,n$ and  $t\in U_\delta$, 
\begin{equation}\label{binodecot}
{N_j-t\theta_j \choose \alpha_j} =\frac{1}{\alpha_j!} \prod_{k=0}^{\alpha_j-1} (N_j-t\theta_j -k).
\end{equation}
It follows that 
\be
\item If $\alpha_j \leq N_j$, then ${N_j-t\theta_j \choose \alpha_j}\big|_{t=0} = {N_j \choose \alpha_j}$ and for  $t\in U_\delta$:
$$\big|{N_j-t\theta_j \choose \alpha_j}\big|\leq \frac{1}{\alpha_j!} \prod_{k=0}^{\alpha_j-1} (N_j+1-k) \leq (N_j+1)!;$$
\item If  $\alpha_j \geq N_j+1$, then ${N_j-t\theta_j \choose \alpha_j}\big|_{t=0} = {N_j \choose \alpha_j}=0$ and 
for  $t\in U_\delta$:
\begin{align*}
\big|{N_j-t\theta_j \choose \alpha_j}\big| \leq \frac{| t \theta_j|}{\alpha_j!} \prod_{k=0}^{N_j-1}(N_j-k+1)  \prod_{k=N_j+1}^{\alpha_j-1}(k-N_j+1)\\
= | t \theta_j|\frac{N_j!}{\alpha_j(\alpha_j-1)\cdots(\alpha_j-N_j+1)}(N_j+1)
\leq (N_j+1) |t \theta_j|.
\end{align*}
\ee
We deduce that for  $t\in U_\delta \setminus\{0\}$,
\begin{eqnarray*}
&&|G_{\Nb, \alphab, \thetab}(t)|\\
&&\ll_{\Nb, \thetab}\frac{\prod_{j\in L(\Nb,\alphab)}|t|}
{\prod_{j\in K(\Nb,\alphab)}|t|\prod_{j\notin K(\Nb,\alphab)}\left|t-
\frac{(N_j+\dots+N_n)+(n+1-j)-(\alpha_j+\dots+\alpha_n)}{\theta_j+\dots+\theta_n}\right|}\\
&&\ll_{\Nb, \thetab} \frac{|t|^{q'-q} }{\prod_{j\not \in K(\Nb, \alphab)}
\left(|t(\theta_j+\dots+\theta_n)- (N_j+\dots+N_n)-(n+1-j)+(\alpha_j+\dots+\alpha_n)|\right)}\\
&&\ll_{\Nb, \thetab} \frac{|t|^{q'-q} }{\prod_{j\not \in K(\Nb, \alphab)}(1/2)}\ll_{\Nb, \thetab}|t|^{q'-q}.
\end{eqnarray*}
It follows that 
$G_{\Nb, \alphab, \thetab}(t)$ is analytic in the whole disk $U_\delta$ and verifies in it the uniform estimate 
$G_{\Nb, \alphab, \thetab}(t) \ll_{\Nb, \thetab}|t|^{q'-q}$.\qed\\
$\bullet$ {\bf Proof of point 3:}  Follows from point 2.\qed\\
$\bullet$ {\bf Proof of point 4:} 
The identity (\ref{binodecot}) implies that if $\alpha_j \geq N_j+1$, then 
$$
{N_j-t\theta_j \choose \alpha_j} \sim_{t\rightarrow 0} \frac{-t \theta_j}{\alpha_j!} \prod_{k=0}^{N_j-1}(N_j-k)  \prod_{k=N_j+1}^{\alpha_j-1}(N_j-k) 
\sim_{t\rightarrow 0}  \frac{t \theta_j (-1)^{\alpha_j-N_j}}{\alpha_j {\alpha_j-1\choose N_j}}.$$
It follows that 
$$G_{\Nb, \alphab, \thetab}(t)\sim_{t\rightarrow 0}\frac{(-1)^{n-q}\left(\prod_{j\in L(\Nb, \alphab)} \frac{(-1)^{\alpha_j -N_j} \theta_j}{\alpha_j {\alpha_j -1 \choose N_j}}\right) 
\prod_{j\not \in  L(\Nb, \alphab)} {N_j \choose \alpha_j}}{\prod_{j\not \in K(\Nb, \alphab)}\left(\frac{(N_j+\dots+N_n)+(n+1-j)-(\alpha_j+\dots+\alpha_n)}{\theta_j+\dots+\theta_n}\right)}.$$
This ends the proof of point 4 and therefore ends the proof of Lemma \ref{keycomputation}. \qed\par

%%%%%%%%%%%%%%%%%%%%%%%%%%%%%%%%%%%%%%%%%%%%%%%%%%%%%%%%%%%%%%%%%%%%%%%%%%%%%%%%%%%%%%%%%%
\section{The key propositions}
%%%%%%%%%%%%%%%%%%%%%%%%%%%%%%%%%%%%%%%%%%%%%%%%%%%%%%%%%%%%%%%%%%%%%%%%%%%%%%%%%%%%%%%%%%%

Now we introduce 
a class of multivariate zeta functions which are slightly more general than that 
considered in Theorem \ref{main1}.    We are working in this slightly more general class because 
it is more suitable for induction arguments.

Let $\q=(q_1,\dots,q_n)\in \N^n$. Set $q=|\q|=q_1+\dots+q_n$.
We will use the notation $\s=(s_{1,1},\dots,s_{1,q_1}, \dots, s_{j,1},\dots,s_{j,q_j},\dots, s_{n,1},\dots,s_{n,q_n})$ for elements of $\C^q$, and denote
$|\s|=s_{1,1}+\cdots+s_{1,q_1}+ \cdots+ s_{j,1}+\cdots+s_{j,q_j}+\cdots+ s_{n,1}+\cdots+s_{n,q_n}$.
Let $\eps \geq 0$ (notice, here, we admit the case $\eps=0$), $\gammab \in \C^n$, and define
\begin{eqnarray}\label{wepsqn}
&&W_\eps (\q,n):=\left\{(\u, \gammab) \in \C^q\times \C^n \mid  \Re(\gamma_j) >\eps 
{\mbox { and }} \Re\left(u_{j,k}+  \gamma_1 \right)>\eps\right.\nonumber\\
&&\qquad\qquad\qquad\qquad\qquad\left. \quad {\rm for\; all}\;\;
 j=1,\dots,n\; {\mbox { and }} \; k=1,\dots, q_j\right\},\\
&& V_{\eps,\q}(\gammab):=\left\{\u \in \C^q\mid \Re\left(u_{j,k}+  \gamma_1 \right)>\eps \quad 
{\rm for\;all}\;\; j=1,\dots,n \;{\mbox { and }} \; k=1,\dots, q_j\right\},\nonumber
\end{eqnarray}
and
$$
\mathcal{D}_{n,\q}:=\left\{\s\in \C^q \;\left|\; \Re\left(\sum_{i=j}^n \sum_{k=1}^{q_i} 
s_{i,k}\right)>n+1-j\quad{\rm for\;all}\;\; j=1,\dots,n\right.\right\}.
$$
For $\s\in \mathcal{D}_{n,\q}$ and $(\u, \gammab) \in  W_0(\q, n)$, define
\begin{equation}\label{EZMaYdef2}
Y_{n,\q}(\s; \u;  \gammab) :=
\int_{[1,\infty)\times [0,\infty)^{n-1}} \prod_{j=1}^n \prod_{k=1}^{q_j} \big(\gamma_1 x_1+\dots+ \gamma_j x_j +u_{j,k}\big)^{-s_{j,k}}~dx_n\dots dx_1
\end{equation}
and
\begin{equation}\label{EZMazetadef2}
Z_{n,\q}(\s; \u;  \gammab) :=
\sum_{m_1\geq 1 \atop m_2, \dots, m_n \geq 0} \frac{1}{\prod_{j=1}^n \prod_{k=1}^{q_j} \big(\gamma_1 m_1+\dots+ \gamma_j m_j +u_{j,k}\big)^{s_{j,k}}}.
\end{equation}
The multiple zeta-function $Z_{n,\q}(\s; \u;  \gammab)$ is absolutely convergent in the region
$\mathcal{D}_{n,\q}$, and in this region
\begin{equation}\label{integral}
Y_{n,\q}(\s; \u;  \gammab) = \int_{[0,1]^n} Z_{n,\q}(\s; \u_{\q}(\b); \gammab)~d\b. 
\end{equation}
Here, $\u_{\q}(\b)\in\mathbb{C}^q$ is given by
$$\u_{\q}(\b)=\left(u_{1,1}(\b),\dots,u_{1,q_1}(\b), \dots, u_{j,1}(\b),\dots,u_{j,q_j}(\b),\dots, u_{n,1}(\b),\dots, u_{n,q_n}(\b)\right),$$ 
where $\b=(b_1,\ldots,b_n)\in[0,1]^n$ and
$u_{j,k}(\b)=u_{j,k} +\sum_{i=1}^j \gamma_i b_i$ for all $j=1,\dots, n$ and all $k=1,\dots, q_j$.

Now we state a proposition, which gives several analytic properties of $Y_{n,\q}(\s; \u;  \gammab)$ and
$Z_{n,\q}(\s; \u;  \gammab)$.

\begin{prop}\label{seriesparameter2}
\be
\item The functions
$\s\mapsto Y_{n,\q}(\s; \u; \gammab)$ and $\s\mapsto Z_{n,\q}(\s; \u; \gammab)$ can be 
meromorphically continued to $\C^q$ and their poles are located in the set 
$$\mathcal P_{n,\q}:=\bigcup_{j=1}^n \bigcup_{k_j \in \N_0}\left\{\s\in \C^q \left| \sum_{i=j}^n \sum_{k=1}^{q_i} s_{i,k}=n+1-j -k_j \right.\right\} .$$
Therefore \eqref{integral} is valid
for all $\s \in \C^n \setminus \mathcal P_{n,\q}$.
\item For  fixed $\omegab \in \C^q$ and  $\thetab \in \C^q$ such that $\sum_{i=j}^n\sum_{k=1}^{q_i}\theta_{i,k} \neq 0$ for all $j=1,\dots, n$, 
there exist $\delta=\delta (\omegab, \thetab)>0$ and $M=M(\omegab, \thetab) >0$ such that 
\be
\item 
$(t, \u, \gammab)\mapsto t^n~Y_{n,\q}(\omegab+ t\thetab; \u; \gammab)$ and  
$(t, \u, \gammab)\mapsto t^n~Z_{n,\q}(\omegab+ t\thetab; \u; \gammab)$
are analytic in the domain 
$\ds U_{\delta}\times  W_0(\q, n)$, where $U_\delta=\{t\in \C;~|t|<\delta\}$;
\item for  $\eps >0$ and  $\gammab \in \mathcal \C^n$ such that $\Re (\gamma_j) > \eps$ 
for all $j=1,\dots, n$,\\
 we have 
$$|t^n~Y_{n,\q}(\omegab+ t\thetab; \u; \gammab)|\ll_{\omegab, \thetab, \gammab, \eps} (1+|\u|)^M$$
and  
$$|t^n~Z_{n,\q}(\omegab+ t\thetab; \u; \gammab)| \ll_{\omegab, \thetab,\gammab, \eps} (1+|\u|)^M$$
uniformly in $(t, \u) \in U_{\delta}\times  V_{\eps,\q}(\gammab)$.
\ee
\ee
\end{prop}
Proposition \ref{seriesparameter2} implies the following key result:
\begin{coro}\label{seriesparameter2cor}
Let $\omegab \in \C^q$ and $\thetab \in \C^q$ be such that $\sum_{i=j}^n\sum_{k=1}^{q_i}\theta_{i,k} \neq 0$ for all $j=1,\dots, n$.
For $(\u, \gammab) \in  W_0(\q, n)$ (see (\ref{wepsqn})), each of the meromorphic functions 
$$t \mapsto Y_{n,\q}(\omegab+ t\thetab; \u; \gammab) {\mbox { and }} 
t\mapsto Z_{n,\q}(\omegab+ t\thetab; \u; \gammab)$$ has at most a pole of order $n$ at $t=0$.
Moreover if we write 
$$ Y_{n,\q}(\omegab+ t\thetab; \u; \gammab) =\sum_{k=0}^n \frac{y_{-k,\q}(\u; \omegab, \thetab; \gammab)}{t^k} +O(t) ~~{\mbox { as }} t\rightarrow 0,$$
 and 
$$ Z_{n,\q}(\omegab+ t\thetab; \u; \gammab) =\sum_{k=0}^n \frac{z_{-k,\q}(\u; \omegab, \thetab; \gammab)}{t^k} +O(t)  ~~{\mbox { as }} t\rightarrow 0.$$
Then, for any $k=0,\dots, n$, 
\be
\item the functions 
$$(\u, \gammab) \mapsto y_{-k,\q}(\u; \omegab, \thetab; \gammab) {\mbox { and }} (\u, \gammab) \mapsto z_{-k,\q}(\u; \omegab, \thetab; \gammab)$$
are analytic in the domain $W_0(\q, n)$ and
\begin{equation}\label{link1}
y_{-k,\q}(\u; \omegab, \thetab; \gammab) = \int_{[0,1]^n} z_{-k,\q}(\u_{\q}(\b); \omegab, \thetab; \gammab)~d\b
\end{equation}
holds in that domain.
\item
There exists  $M=M(\omegab, \thetab)>0$ such that for  $\eps >0$, and $\gammab \in \C^n$ such that $\Re(\gamma_j) > \eps$ for all $j=1,\dots, n$, 
we have uniformly in $\u \in  V_{\eps, \q} (\gammab)$:
$$y_{-k,\q}(\u; \omegab, \thetab; \gammab) \ll_{\omegab, \thetab, \gammab, \eps} (1+|\u|)^M {\mbox { and }} 
z_{-k,\q}(\u; \omegab, \thetab; \gammab) \ll_{\omegab, \thetab,\gammab, \eps} (1+|\u|)^M. $$
\ee
\end{coro}
{\bf Deduction of Corollary \ref{seriesparameter2cor} from Proposition \ref{seriesparameter2}:}\\
The corollary follows from point 2 of Proposition  \ref{seriesparameter2} by applying 
Cauchy's formula which expresses the coefficients of Laurent's expansion of a given one variable meromorphic  function in terms of its integrals on small disks around its singular point.
The identity (\ref{link1}) follows by using in addition the equality \eqref{integral}. \qed
\bigskip

{\bf Proof of Proposition \ref{seriesparameter2}:}\\
The proof  of the proposition for $Y_{n,\q}(\s; \u;  \gammab)$ is similar (and more  easier) than its proof for $Z_{n,\q}(\s; \u;  \gammab)$. 
So we will give here only the proof for $Z_{n,\q}(\s; \u;  \gammab)$.
We will prove the proposition for $Z_{n,\q}(\s; \u;  \gammab)$ by induction on $n$.

$\bullet$ {\bf Proof of Proposition \ref{seriesparameter2} in the case $n=1$:}\\
For  $(\u,\gamma)=((u_1,\ldots,u_q),\gamma) \in W_0(q, 1)$ and  
$\s=(s_1,\dots,s_q)\in \mathcal D_{1, q}$, we have
$$
Z_{1,\q}(\s; \u;  \gamma) =
\sum_{m \geq 1 } \frac{1}{ \prod_{k=1}^{q} \big(\gamma m +u_{k}\big)^{s_{k}}}.$$
Let $K\in \N_0$. Define for $m\geq 1$, the function 
$\ds{ \psi_m(z):=\prod_{k=1}^{q} \big(1  +\frac{u_{k}}{\gamma m} z \big)^{-s_{k}}}$.
Since
$$
\prod_{k=1}^{q} \big(\gamma m +u_{k}\big)^{-s_{k}} 
=\prod_{k=1}^{q}(\gamma m)^{-s_k} \psi_m(1),
$$
applying Taylor's formula with remainder (\cite[(3.4)]{EMT})
to the function $\psi_m(z)$, we obtain that 
for  $m\geq 1$,
\begin{eqnarray*}
&&\prod_{k=1}^{q} \big(\gamma m +u_{k}\big)^{-s_{k}} 
=\sum_{\alphab \in \N_0^q\atop |\alphab|\leq K} {-\s \choose \alphab} \u^{\alphab} \gamma^{-|\s|-|\alphab|}  m^{-|\s|-|\alphab|} \\
&& \quad+(K+1) \sum_{\alphab \in \N_0^q\atop |\alphab|= K+1} {-\s \choose \alphab} \u^{\alphab}\int_0^1 (1-y)^K \prod_{k=1}^{q} \big(\gamma  m +u_{k} y \big)^{-s_{k}-\alpha_k}~dy.
\end{eqnarray*}
It follows that for  $(\u,\gammab) \in W_0(\q, 1)$ and  $\s\in \mathcal D_{1,q}$, 
\begin{eqnarray*}
Z_{1,\q}(\s; \u;  \gammab) &=& \sum_{\alphab \in \N_0^q\atop |\alphab|\leq K} {-\s \choose \alphab} \u^{\alphab} \gamma^{-|\s|-|\alphab|}\zeta (|\s|+|\alphab|)\\
&&+(K+1) \sum_{\alphab \in \N_0^q\atop |\alphab|= K+1} {-\s \choose \alphab} \u^{\alphab}\mathcal R_K(\s; \u; \gamma; \alphab)
\end{eqnarray*}
where 
$$
R_K(\s; \u; \gamma; \alphab) = \sum_{m \geq 1 }\int_0^1 (1-y)^K \prod_{k=1}^{q} \big(\gamma  m +u_{k} y \big)^{-s_{k}-\alpha_k}~dy.
$$
Let $\eps >0$. 
We have uniformly in $m\in \N$, $(\u, \gamma) \in W_\eps (\q, 1)$, $y\in [0,1]$ and 
$k\in \{1,\dots, q\}$:
\begin{eqnarray*}
|\gamma  m +u_{k} y| &\geq& \Re (\gamma) m + \Re(u_k) y =\eps m +\left( \Re(\gamma)-\eps\right) m +\Re( u_k)  y \\
&\geq& \eps m + \left(\Re(\gamma) -\eps+\Re(u_k )\right) y \geq \eps m
\end{eqnarray*}
and
$$|\gamma  m +u_{k} y| \leq  |\gamma|  m + |u_{k}| \leq (|\gamma|+|u_k|) m.$$
The theorem of analyticity under the integral sign implies then that 
$$(\s, \u, \gamma)\mapsto R_K(\s; \u;\gamma;\alphab)$$
is holomorphic in the domain $\{\s \in \C^q \mid \Re(s_1+\dots+s_q)>-K\}\times W_\eps (q, 1)$ and 
the estimate
$$ R_K(\s; \u;\gamma;\alphab) \ll_{s, K, \eps} (1+|\u|+|\gamma|)^{|\s|+K+1}$$
holds there uniformly in $(\u, \gamma) \in  W_\eps (q, 1).$
By using in addition the classical properties of the Riemann zeta function, we deduce that $\s\mapsto Z_{1,\q}(\s; \u;  \gamma)$ has a meromorphic continuation to 
$\{\s \in \C^q \mid \Re(s_1+\dots+s_q)>-K\}$ with poles located in the set $\mathcal P_{1,\q}$ 
and that the point 2 holds for any 
$\omegab \in \C^q$ such that $\Re(\omega_1+\dots +\omega_q )> -K$.\\
By letting $K\rightarrow \infty$ and $\eps \rightarrow 0$, we end the proof of Proposition \ref{seriesparameter2} in the case $n=1$.\qed

$\bullet${\bf Let $n\in \N$ such that $n\geq 2$. We assume that  Proposition \ref{seriesparameter2} holds for $n-1$. We will prove that it remains valid for $n$:}\\
Let $(\u, \gammab) \in W_0(\q, n)$ and $\s\in  \mathcal D_{n,\q}$. 
Fix $m_1\geq 1$ and $m_2,\dots, m_{n-1}\geq 0$.\\
The function $\varphi (x):=\prod_{i=1}^{q_n} \big(\gamma_1 m_1+\dots+ \gamma_{n-1} m_{n-1}+\gamma_n x +u_{n,i}\big)^{-s_{n,i}}$ belongs to $\mathcal C^\infty[0,\infty)$ and 
for all $k\in \N_0$ and all $x\in [0,\infty)$,
$$\varphi^{(k)}(x)=k!  \sum_{\alphab \in \N_0^{q_n}\atop |\alphab|= k} \gamma_n^{|\alphab|} \prod_{i=1}^{q_n}{-s_{n,i} \choose \alpha_i}
\big(\gamma_1 m_1+\dots+ \gamma_{n-1} m_{n-1}+\gamma_n x +u_{n,i}\big)^{-s_{n,i}-\alpha_i}.$$
Let $K\in \N_0$, and let
$\tilde{B}_k$ ($k\geq 0$) be the modified Bernoulli numbers defined by  ${\tilde{B}}_k :=B_k$ 
for all $k\neq 1$ and ${\tilde{B}}_1:=-B_1=\frac{1}{2}$.    (In some references $\tilde{B}_k$ is
written as $B_k$.)
By applying the Euler-Maclaurin formula to the above $\varphi(x)$, we obtain that
\begin{eqnarray}\label{onehand}
&&\sum_{m_n=0}^\infty \prod_{i=1}^{q_n} \big(\sum_{j=1}^n\gamma_j m_j+u_{n,i}\big)^{-s_{n,i}}\nonumber\\
&=& \int_0^\infty \prod_{i=1}^{q_n} \big(\sum_{j=1}^{n-1}\gamma_j m_j+\gamma_n x +u_{n,i}\big)^{-s_{n,i}} ~dx \\
&&+\sum_{\alphab \in \N_0^{q_n}\atop |\alphab|\leq K} 
\frac{(-1)^{|\alphab|} {\tilde{B}}_{|\alphab|+1}\gamma_n^{|\alphab|}}{|\alphab|+1} \prod_{i=1}^{q_n}{-s_{n,i} \choose \alpha_i}
\big(\sum_{j=1}^{n-1}\gamma_j m_j+u_{n,i}\big)^{-s_{n,i}-\alpha_i}\nonumber\\
&&+(-1)^K \gamma_n^{K+1}\sum_{\alphab \in \N_0^{q_n}\atop |\alphab|=K+1} \prod_{i=1}^{q_n}{-s_{n,i} \choose \alpha_i}
\int_0^\infty B_{K+1}(x)~\prod_{i=1}^{q_n} \big(\sum_{j=1}^{n-1}\gamma_j m_j+\gamma_n x +u_{n,i}\big)^{-s_{n,i}-\alpha_i} ~dx,\nonumber
\end{eqnarray}
where $B_{K+1}(x)$ is the $(K+1)$-th periodic Bernoulli polynomial.
On the integrand in the first integral in (\ref{onehand}), again using Taylor's formula with remainder (\cite[(3.4)]{EMT}), we have
\begin{eqnarray}\label{otherhand}
&&\prod_{i=1}^{q_n} \big(\sum_{j=1}^{n-1}\gamma_j m_j+\gamma_n x +u_{n,i}\big)^{-s_{n,i}}  \nonumber\\
&=& \big(\sum_{j=1}^{n-1}\gamma_j m_j+\gamma_n x \big)^{-\sum_{i=1}^{q_n}s_{n,i}} 
\prod_{i=1}^{q_n} \big(1+\frac{u_{n,i}}{\sum_{j=1}^{n-1}\gamma_j m_j+\gamma_n x}\big)^{-s_{n,i}} \\
&=& \sum_{\alphab\in \N_0^{q_n}\atop |\alphab|\leq K} 
\left(\prod_{i=1}^{q_n}{-s_{n,i}\choose \alpha_i} u_{n,i}^{\alpha_i}\right) \big(\sum_{j=1}^{n-1}\gamma_j m_j+\gamma_n x\big)^{-\sum_{i=1}^{q_n}(s_{n,i}+\alpha_i)}\nonumber\\
&&+(K+1) \sum_{\alphab\in \N_0^{q_n}\atop |\alphab|= K+1} 
\left(\prod_{i=1}^{q_n}{-s_{n,i}\choose \alpha_i} u_{n,i}^{\alpha_i}\right)\nonumber\\
&&\qquad\times\int_0^1 (1-y)^K \prod_{i=1}^{q_n} \big(\sum_{j=1}^{n-1}\gamma_j m_j+\gamma_n x + y u_{n,i}\big)^{- s_{n,i}-\alpha_i} ~dy.\nonumber
\end{eqnarray}
Substituting (\ref{onehand}) and (\ref{otherhand}) into (\ref{EZMazetadef2}), and carrying out the first integral, we find that,
for $K\in \N_0$, $(\u, \gammab) \in W_0(\q, n)$ and  $\s\in  \mathcal D_{n,\q}$:
\begin{eqnarray}\label{reductionton-1}
&&Z_{n,\q}(\s; \u;  \gammab) 
= \sum_{m_1\geq 1 \atop m_2, \dots, m_{n-1} \geq 0} \frac{1}{\prod_{j=1}^{n-1} \prod_{k=1}^{q_j} \big(\gamma_1 m_1+\dots+ \gamma_j m_j +u_{j,k}\big)^{s_{j,k}}}\nonumber\\
&&\qquad\times\sum_{m_n=0}^{\infty}\frac{1}{\prod_{i=1}^{q_n}(\gamma_1 m_1+\dots+ \gamma_n m_n +u_{n,i})
^{s_{n,i}}}
\nonumber\\
&&= 
 \sum_{\alphab\in \N_0^{q_n}\atop |\alphab|\leq K} 
\frac{\left(\prod_{i=1}^{q_n}{-s_{n,i}\choose \alpha_i} u_{n,i}^{\alpha_i}\right) }{\gamma_n \left(-1+\sum_{i=1}^{q_n}(s_{n,i}+\alpha_i)\right)}\nonumber\\
&& \times 
\sum_{m_1\geq 1 \atop m_2, \dots, m_{n-1} \geq 0} \frac{1}{\left[\prod_{j=1}^{n-1} \prod_{k=1}^{q_j} \big(\sum_{i=1}^j\gamma_i m_i +u_{j,k}\big)^{s_{j,k}}\right]
\big(\sum_{j=1}^{n-1}\gamma_j m_j\big)^{\sum_{i=1}^{q_n}(s_{n,i}+\alpha_i )-1}}\nonumber\\
&&+\sum_{\alphab \in \N_0^{q_n}\atop |\alphab|\leq K} \frac{(-1)^{|\alphab|} {\tilde{B}}_{|\alphab|+1}\gamma_n^{|\alphab|}}{|\alphab|+1} \prod_{i=1}^{q_n}{-s_{n,i} \choose \alpha_i}\nonumber\\
&& \times 
\sum_{m_1\geq 1 \atop m_2, \dots, m_{n-1} \geq 0} \frac{1}{\left[\prod_{j=1}^{n-1} \prod_{k=1}^{q_j} \big(\sum_{i=1}^j\gamma_i m_i +u_{j,k}\big)^{s_{j,k}}\right]
\big(\prod_{i=1}^{q_n}\big(\sum_{j=1}^{n-1}\gamma_j m_j+u_{n,i}\big)^{s_{n,i}+\alpha_i} \big)}\nonumber\\
&& +\mathcal R_{K,n}^1(\s; \u; \gammab) +\mathcal R_{K,n}^2(\s; \u; \gammab) ,
\end{eqnarray}
where 
\begin{eqnarray*}
&&\mathcal R_{K,n}^1(\s; \u; \gammab) = (K+1) \sum_{\alphab\in \N_0^{q_n}\atop |\alphab|= K+1} 
\left(\prod_{i=1}^{q_n}{-s_{n,i}\choose \alpha_i} u_{n,i}^{\alpha_i}\right) \\
&& \quad\times
\sum_{m_1\geq 1 \atop m_2, \dots, m_{n-1} \geq 0} \frac{\int_0^\infty \int_0^1 (1-y)^K \prod_{i=1}^{q_n} 
\big(\sum_{j=1}^{n-1}\gamma_j m_j+\gamma_n x + y u_{n,i}\big)^{- s_{n,i}-\alpha_i} dy dx}{\prod_{j=1}^{n-1} \prod_{k=1}^{q_j} \big(\gamma_1 m_1+\dots+ \gamma_j m_j +u_{j,k}\big)^{s_{j,k}}},
\end{eqnarray*}
and
\begin{eqnarray*}
&&\mathcal R_{K,n}^2(\s; \u; \gammab) = 
(-1)^K \gamma_n^{K+1}\sum_{\alphab \in \N_0^{q_n}\atop |\alphab|=K+1} \prod_{i=1}^{q_n}{-s_{n,i} \choose \alpha_i}\\
&& \quad\times \sum_{m_1\geq 1 \atop m_2, \dots, m_{n-1} \geq 0} \frac{\int_0^\infty B_{K+1}(x)~\prod_{i=1}^{q_n} \big(\sum_{j=1}^{n-1}\gamma_j m_j+\gamma_n x +u_{n,i}\big)^{-s_{n,i}-\alpha_i} dx}{\prod_{j=1}^{n-1} \prod_{k=1}^{q_j} \big(\gamma_1 m_1+\dots+ \gamma_j m_j +u_{j,k}\big)^{s_{j,k}}}.
\end{eqnarray*}
The formula \eqref{reductionton-1} is the key for the induction process.
In fact, the induction hypothesis implies that the first two terms on the right-hand side of
\eqref{reductionton-1} can be continued meromorphically to the whole space,
and their poles are located in the set $\mathcal P_{n,\q}$.

The remaining task is to evaluate $\mathcal R_{K,n}^1(\s; \u; \gammab)$ and
$\mathcal R_{K,n}^2(\s; \u; \gammab)$.
Define 
$$
\ds \mathcal D_{n,\q}(K):=\{\s\in \C^q \mid \Re\left(\sum_{i=j}^n \sum_{k=1}^{q_i} s_{i,k}\right)>n+1-j -K ~{\rm for\;all}\; j=1,\dots,n\}.
$$
Let $\eps >0$. We have uniformly in $x_1\geq 1$, $x_2,\dots,x_n\geq 0$, $(\u, \gammab) \in W_\eps (\q, n)$, $y\in [0,1]$, $j\in\{1,\dots, n\}$ and $k\in\{1,\dots,q_j\}$:
\begin{eqnarray}\label{minorationsum}
|\gamma_1  x_1+\dots+ \gamma_j x_j +u_{j, k} y| &\geq& \Re (\gamma_1) x_1 + \Re(u_{j, k}) y + \sum_{i=2}^j  \Re (\gamma_i) x_i\nonumber\\
&=&\eps x_1 +\left( \Re(\gamma_1)-\eps\right) x_1 +\Re( u_{j,k})  y + \sum_{i=2}^j  \Re (\gamma_i) x_i \nonumber\\
&\geq& \eps x_1 +\left( \Re(\gamma_1)-\eps +\Re( u_{j,k}) \right) y + \sum_{i=2}^j  \Re (\gamma_i) x_i\nonumber \\
&\geq& \eps x_1  + \sum_{i=2}^j  \Re (\gamma_i) x_i \geq \eps \left(\sum_{i=1}^j x_i\right),
\end{eqnarray}
and
\begin{eqnarray}\label{majorationsum}
|\gamma_1  x_1+\dots+ \gamma_j x_j +u_{j, k} y| &\leq& |\gamma_1|  x_1+\dots+| \gamma_j| x_j +|u_{j, k}| \nonumber\\
&\leq& (1+|\u|+ |\gammab|) \left(x_1+\dots+x_j\right).
\end{eqnarray}
Combining (\ref{minorationsum}) and (\ref{majorationsum}) we see that for any $\eps \in (0,1)$ and any compact subset $H$ of $\C$, 
 we have uniformly in $x_1\geq 1$, $x_2,\dots,x_n\geq 0$, $(\u, \gammab) \in W_\eps (\q, n)$, $y\in [0,1]$, $j\in\{1,\dots, n\}$, $k\in\{1,\dots,q_j\}$ and $s\in H$:
\begin{eqnarray}\label{majorationsumavecs}
\lefteqn{|\left(\gamma_1  x_1+\dots+ \gamma_j x_j +u_{j, k} y\right)^{-s}|}\nonumber\\ 
&=& |\gamma_1  x_1+\dots+ \gamma_j x_j +u_{j, k} y|^{-\Re(s)} e^{\Im (s) \arg \left(\gamma_1  x_1+\dots+ \gamma_j x_j +u_{j, k} y\right)}\nonumber\\
&\leq & |\gamma_1  x_1+\dots+ \gamma_j x_j +u_{j, k} y|^{-\Re(s)} e^{\frac{\pi}{2} |\Im (s)| }\nonumber\\
&\ll_H& |\gamma_1  x_1+\dots+ \gamma_j x_j +u_{j, k} y|^{-\Re(s)} \nonumber\\
&\ll_H& \begin{cases}
 \eps^{-\Re(s)} \left(x_1+\dots+x_j\right)^{-\Re(s)} & \text{if } \Re (s) \geq  0 \\
  (1+|\u|+ |\gammab|)^{-\Re(s)} \left(x_1+\dots+x_j\right)^{-\Re(s)} & \text{if } \Re(s) <  0
  \end{cases} \nonumber \\
  &\ll_{H, \eps}& \begin{cases}
 \left(x_1+\dots+x_j\right)^{-\Re(s)} & \text{if } \Re (s) \geq  0 \\
  (1+|\u|+ |\gammab|)^{|s|} \left(x_1+\dots+x_j\right)^{-\Re(s)} & \text{if } \Re(s) <  0
  \end{cases} \nonumber \\
&\ll_{H, \eps}&  (1+|\u|+ |\gammab|)^{|s|} \left(x_1+\dots+x_j\right)^{-\Re(s)}.
\end{eqnarray}
We deduce that for $K\in \N_0$,  $\eps >0$, $\alphab \in \N_0^{q_n}$ such that $|\alphab|=K+1$, and any compact subset $\mathcal K$ of $\mathcal D_{n,\q}(K)$, 
 we have uniformly in $(\u, \gammab) \in W_\eps (\q, n)$, in $\s \in \mathcal K$ and in $m_1\geq 1$ and $m_2,\dots, m_n\geq 0$:
\begin{eqnarray*}
&& \left| \frac{\int_0^\infty \int_0^1 (1-y)^K \prod_{i=1}^{q_n} 
\big(\sum_{j=1}^{n-1}\gamma_j m_j+\gamma_n x + y u_{n,i}\big)^{- s_{n,i}-\alpha_i} dy dx}{\prod_{j=1}^{n-1} \prod_{k=1}^{q_j} \big(\gamma_1 m_1+\dots+ \gamma_j m_j +u_{j,k}\big)^{s_{j,k}}}\right| \\
&\ll_{K, \mathcal K, \eps}&(1+|\u|+|\gammab|)^{|\s|+|\alphab|}
\frac{\int_0^\infty \int_0^1 (1-y)^K 
\big(\sum_{j=1}^{n-1} m_j+ x \big)^{- \Re\sum_{i=1}^{q_n}( s_{n,i}+\alpha_i)} dy dx}{\prod_{j=1}^{n-1} \prod_{k=1}^{q_j} \big(m_1+\dots+ m_j \big)^{\Re (s_{j,k})}}\\
&\ll_{K, \mathcal K, \eps}&(1+|\u|+|\gammab|)^{|\s|+|\alphab|}
\frac{\int_0^\infty  
\big(\sum_{j=1}^{n-1} m_j+ x \big)^{- \Re\sum_{i=1}^{q_n} (s_{n,i}+\alpha_i)} dx}{\prod_{j=1}^{n-1} \prod_{k=1}^{q_j} \big(m_1+\dots+ m_j \big)^{\Re (s_{j,k})}}\\
&\ll_{K, \mathcal K, \eps}&
\frac{(1+|\u|+|\gammab|)^{|\s|+K+1}}{\left(\prod_{j=1}^{n-1} (m_1+\dots+m_j)^{\Re(\sum_{i=1}^{q_j} s_{j,i})}\right) (m_1+\dots+m_{n-1})^{\Re(\sum_{i=1}^{q_n} s_{j,i}) +K}}.
\end{eqnarray*}
In view of \eqref{domaincv},
the theorem of analyticity  under the integral sign implies then that 
$$(\s, \u, \gammab)\rightarrow \mathcal R_{K,n}^1(\s; \u; \gammab)$$ 
is holomorphic in the domain $\mathcal D_{n,\q}(K)\times W_\eps(\q, n)$ and verifies in it the  estimate 
$$\mathcal R_{K,n}^1(\s; \u; \gammab)\ll_{\s, K, \eps} (1+|\u|+|\gammab|)^{|\s|+K+1} {\mbox { uniformly in }} (\u, \gammab)\in W_\eps (\q, n).$$

A similar argument shows that  
$$(\s, \u, \gammab)\rightarrow \mathcal R_{K,n}^2(\s; \u; \gammab)$$ 
is holomorphic in the domain $\mathcal D_{n,\q}(K)\times W_\eps(\q, n)$ and the  estimate 
$$\mathcal R_{K,n}^2(\s; \u; \gammab)\ll_{\s, K, \eps} (1+|\u|+|\gammab|)^{|\s|+K+1}$$
holds there  uniformly in $(\u, \gammab)\in W_\eps (\q, n)$.

Now we can conclude from (\ref{reductionton-1}) that 
$\s\mapsto Z_{n,\q}(\s; \u; \gammab)$ has the meromorphic continuation to $\mathcal D_{n,q}(K)$ with poles located in the set $\mathcal P_{n,\q}$ and that the point 2 of Proposition \ref{seriesparameter2} holds for any 
$\omegab \in \mathcal D_{n,\q}(K)$.

By letting $K\rightarrow \infty$ and $\eps \rightarrow 0$, we end the proof of   Proposition \ref{seriesparameter2} in the case $n$. This finishes the proof of Proposition \ref{seriesparameter2}. \qed\par
Now we can prove the following necessary result:
\begin{prop}\label{limYninN}
Let $\gammab =(\gamma_1,\dots, \gamma_n) \in \C^n$ be  such that 
$\Re (\gamma_j) >0$ for all $j=1,\dots, n$.
Define
$$V_0(\gammab):=\left\{\u \in \C^n\mid \Re\left(u_{j}+  \gamma_1 \right)>0~{\rm for\;all}\; j=1,\dots,n\right\}$$
For $\u\in V_0(\gammab)$ and  $\s\in \mathcal D_{n}$,
define
\begin{equation}\label{def-Y-n}
Y_n(\s; \u; \gammab):=\int_{(1,\infty)\times (0,\infty)^{n-1}} \prod_{j=1}^n \left(\sum_{i=1}^j \gamma_i x_i+u_j\right)^{-s_j}~dx_n\dots dx_1.
\end{equation}
Let $\thetab=(\theta_1,\dots, \theta_n) \in \C^n$ be such that $\theta_j+\dots+\theta_n \neq 0$ for all $j=1,\dots, n$.
 Then, for any $\Nb =(N_1,\dots, N_n) \in \N_0^n$, the limit 
$$Y_n^{\thetab} (-\Nb; \u;\gammab):=\lim_{t\rightarrow 0}Y_n(-\Nb+t\thetab; \u; \gammab)$$
exists, and we have 
$$Y_n^{\thetab} (-\Nb; \u; \gammab) 
=\sum_{I\subset \{1,\dots, n\}} \sum_{\alphab \in \mathcal J (I, \Nb)} A(\Nb, I, \alphab, \thetab, \gammab) ~\u^{\alphab}$$
where 
\begin{eqnarray}\label{groscoefs}
A(\Nb, I, \alphab, \thetab, \gammab) 
&=&
  \frac{ (-1)^{n-|I|+\sum_{j\in L(\Nb, \alphab)} (\alpha_j -N_j)} \prod_{j \not \in L(\Nb, \alphab)} {N_j \choose \alpha_j}}{\prod_{j\in L(\Nb, \alphab)}
 \alpha_j {\alpha_j -1\choose N_j}
 \prod_{j\not \in I} 
 \left(n+1-j+\sum_{i=j}^n N_i-\sum_{i=j}^n\alpha_i\right)}\nonumber\\
 &&\qquad \times \left(\gamma_1^{|\Nb|-|\alphab|+n}\prod_{j=1}^n \gamma_j^{-1}\right)~\left(\frac{ \prod_{j  \in L(\Nb, \alphab)} \theta_j}{\prod_{j \in I}(\theta_j+\dots+\theta_n)}\right).
\end{eqnarray}
\end{prop}
{\bf Proof of Proposition \ref{limYninN}:}\\
First we recall from Proposition \ref{seriesparameter2} that 
$Y_n(\s; \u; \gammab)$ has a meromorphic continuation to the whole complex space $\C^n$ and its poles are located in the set 
$$
\mathcal P_n:=\bigcup_{j=1}^n\bigcup_{k_j \in \N_0} \{\s=(s_1,\dots, s_n) \in \C^n\mid s_j+\dots +s_n =n+1-j -k_j\}.
$$
Define 
\begin{equation}\label{condition}
V_1(\gamma_1):=\{\u \in \C^n \mid |u_j|< \Re (\gamma_1)~~{\rm for\;all}\; j=1,\dots,n\}.
\end{equation}
Let $\s=(s_1,\dots, s_n) \in \mathcal D_n$ and we first assume that $\u \in  V_1( \gamma_1)$.
We have uniformly in $\x=(x_1,\dots,x_n) \in [1,\infty)\times [0,\infty)^{n-1}$:
$$\left|\frac{u_j}{\sum_{i=1}^j \gamma_i x_i}\right|\leq \frac{|u_j|}{\sum_{i=1}^j \Re(\gamma_i) x_i}
\leq \frac{|u_j|}{ \Re(\gamma_1) }<1.$$
Therefore
$$\prod_{j=1}^n \left(\sum_{i=1}^j \gamma_i x_i +u_j\right)^{-s_j}= \sum_{\alphab \in \N_0^n} {-\s\choose \alphab} 
\u^{\alphab} \prod_{j=1}^n \left(\sum_{i=1}^j \gamma_i x_i \right)^{-s_j-\alpha_j},$$
where the right-hand side converges uniformly in $\x=(x_1,\dots,x_n) \in [1,\infty)\times [0,\infty)^{n-1}$,
This implies that for any $\s=(s_1,\dots, s_n) \in \mathcal D_n$,
$$Y_n(\s; \u; \gammab) = \sum_{\alphab \in \N_0^n} {-\s\choose \alphab} 
\u^{\alphab} Y_n(\s+\alphab; \gammab),$$
where
$Y_n(\s; \gammab)$ is defined by \eqref{ynsgammaadef}.
Applying Lemma \ref{integralexpression} we obtain that for any $\s=(s_1,\dots, s_n) \in \mathcal D_n$,
\begin{equation}\label{ynsagammabaexp}
Y_n(\s; \u; \gammab) = \sum_{\alphab \in \N_0^n} 
\frac{\gamma_1^{-|\s|-|\alphab|+n} {-\s\choose \alphab} 
\u^{\alphab}}{(\gamma_1\dots \gamma_n) \prod_{j=1}^n (s_j+\dots+s_n +\alpha_j+\dots+\alpha_n+j-n-1)}.
\end{equation}
Moreover, since $\u \in  V_1( \gamma_1)$, the right-hand side of (\ref{ynsagammabaexp}) is uniformly convergent in any compact subset of  $\C^n \setminus \mathcal P_n$.
It follows that the meromorphic continuation of $Y_n(\s; \u; \gammab)$  is given by (\ref{ynsagammabaexp}) for any $\s \in \C^n\setminus \mathcal P_n$. \par
Let $\thetab=(\theta_1,\dots, \theta_n) \in \C^n$ such that 
$\sum_{i=j}^n\theta_i \neq 0$ for all $j$ and $\Nb =(N_1,\dots, N_n) \in \N_0^n$.
Set $\delta :=\frac{1}{2} \min\left\{(1+|\theta_j|)^{-1}, |\theta_j+\dots+\theta_n|^{-1};~j=1,\dots, n\right\} \in (0,1/2)$ and $U_\delta =\{t\in \C;~|t|<\delta\}$.
From (\ref{ynsagammabaexp}) we obtain that 
$$Y_n(-\Nb+t\thetab; \u; \gammab) = \sum_{\alphab \in \N_0^n} 
\frac{\gamma_1^{|\Nb|+n-|\alphab| -t|\thetab|} 
\u^{\alphab}}{(\gamma_1\dots \gamma_n) \prod_{j=1}^n (\theta_j+\dots+\theta_n)} G_{\Nb, \alphab, \thetab}(t)$$
for any $t\in U_{\delta}\setminus\{0\}$, 
where $G_{\Nb, \alphab, \thetab}(t)$ is defined by \eqref{def-G}.
By using point 2 of Lemma \ref{keycomputation}, it follows from  Lebesgue's dominated convergence theorem that 
$\ds Y_n^{\thetab} (-\Nb; \u;\gammab):=\lim_{t\rightarrow 0}Y_n(-\Nb+t\thetab; \u; \gammab)$ exists and that 
\begin{equation}\label{formulaaa}
Y_n^{\thetab} (-\Nb; \u;\gammab)= \sum_{\alphab \in \N_0^n} 
\frac{\gamma_1^{|\Nb|+n-|\alphab| } 
\u^{\alphab}}{(\gamma_1\dots \gamma_n) \prod_{j=1}^n (\theta_j+\dots+\theta_n)} ~G_{\Nb, \alphab, \thetab}(0),
\end{equation}
where $G_{\Nb, \alphab, \thetab}(0)$ is defined in Lemma \ref{keycomputation}. 
Moreover, points 3 and 4(b) of Lemma \ref{keycomputation} imply that $G_{\Nb, \alphab, \thetab}(0)=0$ if $\alphab \not\in \{0, |\Nb|+n\}^n$. 
It follows that the sum on the right-hand side of (\ref{formulaaa}) is finite. 

Therefore by using the expression of $G_{\Nb, \alphab, \thetab}(0)$ given by Lemma \ref{keycomputation} and by arranging the terms we obtain that 
\begin{eqnarray}\label{crucialequality}
Y_n^{\thetab} (-\Nb; \u;\gammab)&=& \sum_{\alphab \in \N_0^n\atop |K(\Nb, \alphab)|=|L(\Nb, \alphab)|} 
\frac{\gamma_1^{|\Nb|+n-|\alphab| } 
\u^{\alphab}}{(\gamma_1\dots \gamma_n) \prod_{j=1}^n (\theta_j+\dots+\theta_n)}\nonumber\\
&& \times \frac{(-1)^{n-|K(\Nb, \alphab)|}\left(\prod_{j\in L(\Nb, \alphab)} \frac{(-1)^{\alpha_j -N_j} \theta_j}{\alpha_j {\alpha_j -1 \choose N_j}}\right) 
\prod_{j\not \in  L(\Nb, \alphab)} {N_j \choose \alpha_j}}{\prod_{j\not \in K(\Nb, \alphab)}\left(\frac{(N_j+\dots+N_n)+(n+1-j)-(\alpha_j+\dots+\alpha_n)}{\theta_j+\dots+\theta_n}\right)}\nonumber\\
&=& \sum_{I\subset \{1,\dots, n\}} \sum_{\alphab \in \mathcal J(I, \Nb)} 
A(\Nb, I, \alphab, \thetab, \gammab)~\u^{\alphab},
\end{eqnarray}
where the coefficients $A(\Nb, I, \alphab, \thetab, \gammab)$ are defined by (\ref{groscoefs}).\par
Fix $\Nb\in \N_0^n$. 
We will now extend the region of $\u$ for which the proposition holds.
Denote the last member of \eqref{crucialequality} by $\psi(\u)$. Since $\psi(\u)$ is polynomial in $\u$, it is  analytic on the set $V_0(\gammab)$.
Moreover, Corollary \ref{seriesparameter2cor} implies that for any $\u \in V_0(\gammab)$
$$Y_n^{\thetab} (-\Nb+t\thetab; \u;\gammab)=\sum_{k=0}^n \frac{y_{-k}(\u; -\Nb, \thetab, \gammab)}{t^k} +O(t) {\mbox { as }} t\rightarrow 0,$$
where for any $k=0,\dots, n$, $\u \mapsto y_{-k}(\u; -\Nb, \thetab, \gammab)$ is analytic in the domain $V_0(\gammab)$.\\
On the other hand, (\ref{crucialequality}) implies  that for any $\u \in V_1( \gamma_1)$,
\begin{equation}\label{crucialequalitybis}
y_0(\u; -\Nb, \thetab, \gammab) = \psi (\u) {\mbox { and }} y_{-k}(\u; -\Nb, \thetab, \gammab) =0 ~{\rm for \;all}\; k=1,\dots, n.
\end{equation}
Since  $V_1(\gamma_1)$ is a non-empty open subset of the convex  (and hence connected) open set $V_0(\gammab)$, 
it follows then by analytic continuation that (\ref{crucialequalitybis}) holds  for any  $\u\in V_0( \gammab)$. 
This  ends the proof of  Proposition \ref{limYninN}.\qed 

%%%%%%%%%%%%%%%%%%%%%%%%%%%%%%%%%%%%%%%%%%%%%%%%%%%%%%%%%%%%%%%%%%%%%%%%%%%%%%%%%%%%%%%
\section{An extension of Raabe's lemma}
%%%%%%%%%%%%%%%%%%%%%%%%%%%%%%%%%%%%%%%%%%%%%%%%%%%%%%%%%%%%%%%%%%%%%%%%%%%%%%%%%%%%%%%

Define for any $\delta \in \R$, 
$$\mathcal H_n(\delta):=\{\z=(z_1,\dots, z_n)\in \C^n \mid \Re (z_i) > \delta ~~{\rm for\;all}\; i=1,\dots,n\}.$$

\begin{lemma}\label{Raabeextension} {\bf (An extension of Raabe's lemma)}\\
Let $\delta >0$. Let $g :\mathcal H_n(-\delta) \rightarrow \C$ be an analytic function  in $\mathcal H_n(-\delta)$  
such that there exists two constants $K>0$ and $c \in (0, \pi)$ such that 
\begin{equation}\label{majorationgz}
|g(\z)| \leq K e^{c(|z_1|+\dots+|z_n|)} \quad \forall \z \in  \mathcal{H}_n(-\delta).
\end{equation}
Define  for any $\x=(x_1,\ldots,x_n) \in \mathcal H_n(-\delta)$,
\begin{equation}\label{RaabeTransform}
f(\x)=\int_{[0,1]^n} g(\x +\y) ~d\y.
\end{equation}
Assume that $f$ is a polynomial of degree at most $d$.
Then,
$g$ is also a polynomial of degree at most $d$.
Moreover, if we write $f(\x)=\sum_{\alphab} a_{\alphab} \x^{\alphab}=\sum_{\alphab} a_{\alphab} \prod_{i=1}^n x_i^{\alpha_i}$ (with $\alphab=(\alpha_1,\ldots,\alpha_n)$),
then 
\begin{equation}\label{raabetransfert}
g(\x) =\sum_{\alphab} a_{\alphab} \prod_{i=1}^n B_{\alpha_i}(x_i),
\end{equation}
where the $B_k(x)$ are the classical Bernoulli polynomials. 
\end{lemma}
{\bf Remark:} 
Raabe's transform (\ref{RaabeTransform}) is an important operator which makes it possible to derive several properties of a Dirichlet series from its associated Dirichlet  integral.  
For the history of Raabe's formula, see E. Friedman and S. Ruijsenaars \cite[p.367]{FR}.
E. Friedman and A. Pereira \cite{fredman} proved this lemma under the assumption that
both $f$ and $g$ are polynomials.    For our aim in the present paper, we only assume in Lemma  \ref{Raabeextension} that $g$ is an analytic function in a suitable domain satisfying the estimate
(\ref{majorationgz}) which is necessary for Carlson's theorem that we used in our proof.
A question that deserves more investigation is to find the optimal constant $c$ in (\ref{majorationgz}) for which Lemma \ref{Raabeextension} remains valid.\par
\bigskip

{\bf Proof of Lemma \ref{Raabeextension}:}\\
We will proceed by induction on n:\\
$\bullet$ {\bf The case $n=1$:}\\
The theorem of differentiation under the integral sign implies that 
for any $x \in \mathcal H_1(-\delta)$,
$$0= f^{(d+1)}(x)=\int_{[0,1]} g^{(d+1)}(x +y) ~dy=g^{(d)}(x+1)-g^{(d)}(x).$$
It follows that 
$$g^{(d)}(k)=g^{(d)}(0) \quad {\rm for\;all}\; k\in \N_0.$$
Let $z \in \C$ such that $\Re (z)\geq 0$.    The Cauchy formula and (\ref{majorationgz}) imply that 
$$|g^{(d)}(z)| =\left|\frac{d!}{2\pi i} \int_{|t-z|=\delta/2} \frac{g(t)}{(t-z)^{d+1}}~dt\right|\leq K' e^{c|z|},$$
where $K'=K ~d! \left(\frac{\delta}{2}\right)^{-d} e^{c \delta/2} >0$.

Then it follows from Carlson's classical theorem (F. Carlson \cite{carlson}; see 5.81 in page 186 of \cite{Titchmarsh}) that 
$$g^{(d)} (z) =g^{(d)} (0) \quad {\rm for\;all}\; z\in \C {\mbox { such that }} \Re (z) \geq 0.$$
Thus, $g$ is a polynomial of degree at most $d$.

Now since we know that $f$ and $g$ are both polynomials,  (\ref{raabetransfert}) is a consequence of the  Lemma of Friedman and Pereira (see Lemma 2.4 of \cite{fredman}) of Raabe type.
This ends the proof of Lemma \ref{Raabeextension} in the case $n=1$.\qed

$\bullet$ {\bf Let $n\in \N$. Assume that Lemma \ref{Raabeextension} is true for function in $n-1$ variables, we will prove that it remains valid for function in $n$ variables:}

Let $\delta >0$ and let $g :\mathcal H_n(-\delta) \rightarrow \C$ be an analytic function  satisfying the assumptions of Lemma \ref{Raabeextension}.
Let $\betab \in \N_0^n$ such that $|\betab| >d$.
The Cauchy formula and (\ref{majorationgz}) imply that there exists 
 $K'>0$ and $c \in (0, \pi)$ such that 
\begin{equation}\label{majorationbetagz}
|\partial^{\betab}g(\z)| \leq K' e^{c(|z_1|+\dots+|z_n|)} \quad {\rm for\;all}\; \z \in  
\mathcal{H}_n(-\delta/2).
\end{equation}
Fix $\z'=(z_1,\dots, z_{n-1}) \in \mathcal H_{n-1}(-\delta/2)$ and 
define $h :\mathcal H_1(-\delta/2) \rightarrow \C$ by 
$$h(z_n):= \int_{[0,1]^{n-1}} \partial^{\betab} g(z_1+a_1,\dots, z_{n-1}+a_{n-1}, z_n) ~da_1\dots da_{n-1}.$$
It is easy to see  that $h$ is analytic in $\mathcal H_1(-\delta/2)$ and that (\ref{majorationbetagz}) implies that 
$$|h(z_n)|\leq  K'(\z') ~e^{c|z_n|} \quad {\rm for\;all}\; z_n \in \mathcal H_1(-\delta/2),$$
where $K' (\z') =K' \left(\frac{e^c -1}{c}\right)^{n-1} e^{c(|z_1|+\dots+|z_{n-1}|)}>0$.

On the other hand, since $|\betab|>d$, 
we have for any $z_n \in \mathcal H_1(-\delta/2)$, 
\begin{eqnarray*}
\int_{[0,1]} h(z_n+a_n) ~da_n &=& \int_{[0,1]^{n}} \partial^{\betab} g(z_1+a_1,\dots, z_{n}+a_{n}) ~da_1\dots da_{n} \\
&=& \partial^{\betab} 
\left(\int_{[0,1]^{n}} g(z_1+a_1,\dots, z_{n}+a_{n}) ~da_1\dots da_{n}\right) \\
&=& \partial^{\betab}f(\z)=0.
\end{eqnarray*}
The case $n=1$ implies then that for any $z_n \in \mathcal H_1(-\delta/2)$, $h(z_n)=0$. 
As a conclusion we proved that 
\begin{equation}\label{enfin1}
 \int_{[0,1]^{n-1}} \partial^{\betab} g(z_1+a_1,\dots, z_{n-1}+a_{n-1}, z_n) ~da_1\dots da_{n-1} =0 
\end{equation}
for all $\z=(z_1,\dots, z_n) \in \mathcal H_{n}(-\delta/2)$.

Now fix $z_n \in \mathcal H_1(-\delta/2)$ and define $\ell :  \mathcal H_{n-1}(-\delta/2) \rightarrow \C$ by 
$$\ell (z_1,\dots, z_{n-1}) = \partial^{\betab}g(z_1,\dots, z_{n-1}, z_n).$$
It is easy to see that $\ell$ is analytic in  $\mathcal H_{n-1}(-\delta/2)$ and that 
 (\ref{majorationbetagz}) implies that 
$$
|\ell (\z')| \leq K''(z_n)~e^{c(|z_1|+\dots+|z_{n-1}|)} \quad {\rm for\;all}\; \z'=(z_1,\dots, z_{n-1}) \in  \mathcal{H}_{n-1}(-\delta/2),
$$
where $ K''(z_n) = K' e^{c|z_n|} >0$.

It follows then from our induction hypothesis and (\ref{enfin1}) that 
$$\ell (z_1,\dots, z_{n-1}) =0 \quad {\rm for\;all}\; \z'=(z_1,\dots, z_{n-1}) \in  
\mathcal{H}_{n-1}(-\delta/2)$$
and hence that for any  $\betab \in \N_0^n$ with $|\betab| >d$ we have
$$\partial^{\betab}g(z_1,\dots, z_{n-1}, z_n) =0 \quad  {\rm for\;all}\; \z=(z_1,\dots, z_{n}) \in  \mathcal{H}_{n}(-\delta/2).$$
It follows that $g$ is a polynomial of degree at most $d$. 
Now since we know that both $f$ and $g$ are polynomials,  (\ref{raabetransfert}) is again a 
consequence of Raabe's Lemma of Friedman and Pereira.
This ends the induction argument and the  proof of Lemma \ref{Raabeextension}.\qed \par

We end this section with the following useful lemma.   This lemma is maybe not new. 
But we give a proof of it in order to be self-contained.
For $\deltab =(\delta_1,\dots, \delta_n) \in \R^n$, define
$$\mathcal H_n(\deltab):=\{\z=(z_1,\ldots,z_n) \in \C^n \mid \Re (z_j) >\delta_j ~~{\rm for\;all}\; j=1,\dots, n\}.$$

\begin{lemma}\label{unicitylemma}
Let $\deltab =(\delta_1,\dots, \delta_n) \in \R^n$ and $\mub=(\mu_1,\dots, \mu_n) \in \R^n$ such that $\mu_j \geq \delta_j$ for all $j=1,\dots, n$.
Let $f: \mathcal H_n(\deltab)\rightarrow \C$ be an analytic function.
Assume that $\ds f(\x)=0$ for all $\x \in \prod_{j=1}^n (\mu_j, \infty)$.
Then 
$f(\z) =0$ for all $\z \in \mathcal H_n(\deltab).$
\end{lemma}
{\bf Proof of Lemma \ref{unicitylemma}:}\\
We will prove the lemma by induction on $n$.

If $n=1$ the lemma is clear.
Let $n\geq 2$. Assume that Lemma \ref{unicitylemma} is true for functions of $n-1$ variables. We will prove that it remains true for functions of $n$ variables.

Fix $x_1,\dots, x_{n-1} \in \R$ such that $x_i >\mu_i$ for all $i=1,\dots, n-1$.
Define the function $F: \mathcal{H}_1(\mu_n) \rightarrow \C$ by $F(z)=f(x_1,\dots, x_{n-1}, z)$.
It follows from our assumptions that $F$ is a one variable analytic function  in the domain $\mathcal{H}_1(\mu_n) $ and that $F(x)=0$ for all $x\in (\mu_n, \infty)$.
We deduce then that $F(z)=0$ for all $z \in \mathcal{H}_1(\mu_n) $.
That is, we have
\begin{eqnarray}\label{traceroute}
f(x_1,\dots, x_{n-1}, z_n) =0\;
 {\rm for \;all}\; (x_1,\dots, x_{n-1}, z_n) \in \left(\prod_{j=1}^{n-1} (\mu_j, \infty)\right) \times \mathcal{H}_1(\mu_n) .
\end{eqnarray}

Now fix $z_n \in \C$ such that $\Re (z_n) >\mu_n$.
Let $\mub'=(\mu_1,\ldots,\mu_{n-1})$ and define 
$g:\mathcal H_{n-1}(\mub')\rightarrow \mathbb{C}$
 by 
$g(z_1,\dots, z_{n-1})= f(z_1,\dots, z_{n-1}, z_n)$.
Then $g$ is analytic in $\mathcal H_{n-1}(\mub')$ and (\ref{traceroute}) implies that $g(x_1,\dots, x_{n-1}) =0$ for all $(x_1,\dots, x_{n-1}) \in \prod_{j=1}^{n-1} (\mu_j, \infty)$.
The induction hypothesis implies then that 
$$g(z_1,\dots, z_{n-1})=0 \quad {\rm for\;all}\; (z_1,\dots, z_{n-1}) \in \mathcal H_{n-1}(\mub').$$
We deduce that 
$$f(z_1,\dots, z_n)= 0 \quad {\rm for\;all}\; \z \in \mathcal H_n(\mub).$$
This ends the proof of  Lemma \ref{unicitylemma} since $\mathcal H_n(\mub)$ is a non-empty open subset of the domain $\mathcal H_n(\deltab)$.\qed

%%%%%%%%%%%%%%%%%%%%%%%%%%%%%%%%%%%%%%%%%%%%%%%%%%%%%%%%%%%%%%%%%%%%%%%%%%%%%%%%%%%%%%%%%
\section{Completion of the proof of Theorem \ref{main1}}
%%%%%%%%%%%%%%%%%%%%%%%%%%%%%%%%%%%%%%%%%%%%%%%%%%%%%%%%%%%%%%%%%%%%%%%%%%%%%%%%%%%%%%%%%
Fix $\Nb =(N_1,\dots, N_n)\in \N_0^n$ and $\thetab=(\theta_1,\dots, \theta_n)\in \C^n$.
Assume that 
$$ \theta_j+\dots+\theta_n \neq 0 \quad {\rm for\;all}\; j=1,\dots, n.$$
Set 
\begin{equation*}
W:=\left\{(\u, \gammab) \in \C^n\times \C^n\mid \Re (\gamma_j) >0 {\mbox { and }} \Re\left(u_{j}+  \gamma_1 \right)>0~{\rm for\;all}\; j=1,\dots,n\right\}.
\end{equation*}
For $(\u,\gammab) \in W$ and  $\s\in \mathcal{D}_n$, we consider
$Y_n(\s; \u; \gammab)$ defined by \eqref{def-Y-n} 
and
\begin{equation}\label{def-Z-n}
Z_n(\s; \u; \gammab) :=
\sum_{m_1\geq 1 \atop m_2, \dots, m_n \geq 0} \frac{1}{\prod_{j=1}^n (\gamma_1 m_1+\dots+ \gamma_j m_j +u_j)^{s_j}}.
\end{equation}
We know from Corollary  \ref{seriesparameter2cor} that for any $(\u, \gammab) \in  W$, the  functions 
$$t \mapsto Y_n(-\Nb+ t\thetab; \u; \gammab) {\mbox { and }} t\mapsto Z_n(-\Nb+ t\thetab; \u; \gammab)$$ are meromorphic and have at most a pole of order $n$ at $t=0$.
Write 
$$ Y_n(-\Nb+ t\thetab; \u; \gammab) =\sum_{k=0}^n \frac{y_{-k}(\u; -\Nb, \thetab; \gammab)}{t^k} +O(t) ~~{\mbox { as }} t\rightarrow 0,$$
 and 
$$ Z_n(-\Nb+ t\thetab; \u; \gammab) =\sum_{k=0}^n \frac{z_{-k}(\u; -\Nb, \thetab; \gammab)}{t^k} +O(t)  ~~{\mbox { as }} t\rightarrow 0.$$
Corollary \ref{seriesparameter2cor} implies then that for any $k=0,\dots, n$, the functions 
\begin{equation}\label{content}
(\u, \gammab) \mapsto y_{-k}(\u; -\Nb, \thetab; \gammab) {\mbox { and }} (\u, \gammab) \mapsto z_{-k}(\u; -\Nb, \thetab; \gammab)
\end{equation}
are analytic in the domain $W$ and 
\begin{equation}\label{link1bis}
y_{-k}(\u; -\Nb, \thetab; \gammab) = \int_{[0,1]^n} z_{-k}(\u(\b); -\Nb, \thetab; \gammab)~d\b
\end{equation}
holds in that domain, where 
$\u(\b)=\left(u_{1}(\b),\dots,u_{n}(\b)\right) \in \C^n$ with $\b=(b_1,\ldots,b_n)$, 
$u_{j}(\b)=u_{j} +\sum_{i=1}^j \gamma_i b_i$ for all $j$.

For $\gammab \in \C^n$ such that $\Re (\gamma_j) >0$ for all $j=1,\dots, n$, define 
$$\mathcal V (\gammab):=\left\{\u \in \C^n\;\left|\; \Re\left(u_{j}+  \gamma_1 \right)> \Re (\sum_{i=1}^j \gamma_i) +1\quad {\rm for\;all}\; j=1,\dots,n\right.\right\}.$$
Temporarily we assume that
$\gammab \in (1,\infty)^n$ and $\u \in \mathcal V ( \gammab)$. 
It is easy to see that for all $\a \in \mathcal H_n(-1)$, and all $j=1,\dots,n$,
$$\Re\left(\u_j(\a) +\gamma_1\right)= \Re(u_j +\gamma_1)+ \sum_{i=1}^j \gamma_i \Re(a_i) >\Re(u_j +\gamma_1)- \sum_{i=1}^j \gamma_i >1,$$
that is for all $\a \in \mathcal H_n(-1)$, 
$$\u(\a) \in V_1(\gammab)=\left\{\z \in \C^n \mid \Re\left(z_{j}+  \gamma_1 \right)>1 \quad {\rm for\;all} \;j=1,\dots,n \right\}.$$
Define for  $k=0,\dots, n$ and  $\a \in  \mathcal H_n(-1)$,
$$f_k(\a):=y_{-k}(\u(\a); -\Nb, \thetab; \gammab) \quad {\mbox { and }} \quad g_k(\a):= z_{-k}(\u(\a); -\Nb, \thetab; \gammab).$$
Corollary  \ref{seriesparameter2cor} implies then  that for $k=0,\dots,n$ the following 
three points hold:
\be
\item   $f_k$ and $g_k$ are analytic functions  in $\mathcal H_n(-1)$;
\item $ \ds f_k(\x)=\int_{[0,1]^n} g_k(\x +\y) ~d\y \quad {\rm for\;all}\; \x \in \mathcal H_n(-1);$
\item there exists a constant $M=M(\Nb, \thetab)>0$ such that,
 uniformly in $\a \in \mathcal H_n(-1)$, we have
$g_k(\a) \ll_{\Nb, \thetab, \u, \gammab} (1+|\a|)^M$.
\ee 
On the other hand, Proposition \ref{limYninN} implies that 
for all $\a \in \mathcal H_n(-1)$,
$$f_k(\a) =0 \quad {\rm for\;all}\; k=1,\dots, n$$
and 
\begin{eqnarray*}
f_0(\a) =
\sum_{I\subset \{1,\dots, n\}} \sum_{\alphab \in \mathcal J (I, \Nb)} 
A(\Nb, I,\alphab,\thetab;\gammab)
\prod_{j=1}^n\left(u_{j} +\sum_{i=1}^j \gamma_i a_i\right)^{\alpha_j},
\end{eqnarray*}
where the coefficients $A(\Nb, I, \alphab,\thetab;\gammab)$ are defined by (\ref{groscoefs}).\\
In particular, $f_0(\a)$ is a polynomial in $\a$.

We deduce then from Lemma \ref{Raabeextension} that 
for all $\a \in \mathcal H_n(-1)$, 
$$g_k(\a) =0 \quad {\rm for\;all}\; k=1,\dots,n$$
and 
\begin{eqnarray*}
g_0(\a) &=&
\sum_{I\subset \{1,\dots, n\}} \sum_{\alphab \in \mathcal J (I, \Nb)} \sum_{\k \in \N_0^n,\atop  |\k|\leq |\alphab|} \widetilde{c}_n(\u;\alphab, \k) A(\Nb, I, \alphab,\thetab;\gammab)
 \prod_{j=1}^n B_{k_j}(a_j),
\end{eqnarray*}
where the polynomials $\widetilde{c}_n(\u;\alphab, \k)$  are defined by
\begin{equation}\label{def-widetilde-c}
\prod_{j=1}^n (\sum_{i=1}^j \gamma_iX_i +u_j)^{\alpha_j}= \sum_{\k \in \N_0^n,\atop  |\k|\leq |\alphab|} \widetilde{c}_n(\u;\alphab, \k) ~\X^\k = \sum_{\k \in \N_0^n,\atop  |\k|\leq |\alphab|}
 \widetilde{c}_n(\u;\alphab, \k) ~X_1^{k_1} \dots X_n^{k_n}.
 \end{equation}
 By taking $\a=0$, we obtain that,
 for all $\gammab \in (1, \infty)^n$ and all $\u \in \mathcal V ( \gammab)$:
 \begin{equation}\label{concl1}
  z_{-k}(\u; -\Nb, \thetab; \gammab) = 0  \quad {\rm for\;all}\; k=1,\dots, n
  \end{equation}
and 
\begin{equation}\label{concl2}
z_{0}(\u; -\Nb, \thetab; \gammab) =
\sum_{I\subset \{1,\dots, n\}} \sum_{\alphab \in \mathcal J (I, \Nb)} \sum_{\k \in \N_0^n,\atop  |\k|\leq |\alphab|} \widetilde{c}_n(\u;\alphab, \k)A(\Nb, I, \alphab,\thetab;\gammab) 
 \prod_{j=1}^n B_{k_j}.
\end{equation}
Since $\widetilde{c}_n(\u;\alphab, \k)=c_n(\u;\alphab, \k)\gamma_1^{k_1}\cdots\gamma_n^{k_n}$
(where $c_n(\u;\alphab, \k)$ is defined by \eqref{benoulliexpension}),
the right-hand side of \eqref{concl2} coincides with the right-hand side of
\eqref{thm-1-formula}.

Moreover, for any fixed $\gammab \in (1,\infty)^n$, $\mathcal V (\gammab)$ is a non-empty open subset of the domain $V_0 (\gammab)$ and we know that for all 
$k=0,\dots, n$, $\u\mapsto z_{-k}(\u; -\Nb, \thetab; \gammab)$ is analytic in $V_0 (\gammab)$. 
It follows then by analytic continuation that for any $\gammab \in (1,\infty)^n$ the identities (\ref{concl1}) and (\ref{concl2}) hold for any $\u$ in the whole domain $V_0 (\gammab)$. 

Now fix $\u \in \C^n$ and set $\ds \eta(\u) :=\max \left\{0,  -\Re(u_1),\dots, -\Re (u_n)\right\}$.
Define $$\mathcal G (\u) :=\{ \gammab \in \C^n \mid  \Re (\gamma_1) > \eta(\u) {\mbox { and }} \Re(\gamma_j) >0 ~~{\rm for\;all}\; j=2,\dots, n\}.$$
From the definition of $V_0(\gammab)$ and $W$, it is easy to see that $\{\u\}\times \mathcal G (\u) \subset W$.
It follows then from (\ref{content}) that 
$$\gammab \mapsto z_{-k}(\u; -\Nb, \thetab; \gammab)
{\mbox { is analytic in the domain }} \mathcal G (\u).$$
We already know from the above that the identities (\ref{concl1}) and (\ref{concl2}) hold for $\gammab \in (1,\infty)^n \cap \mathcal G (\u)$.
 Lemma \ref{unicitylemma} implies then that for any $\u \in \C^n$  the identities (\ref{concl1}) and (\ref{concl2}) hold for any $\gammab$ in the whole domain $\mathcal G (\u)$.
This ends the proof of Theorem  \ref{main1}.\qed

\medskip
{\bf Acknowledgments:} The authors wishe to express their thanks to an anonymous referee for 
careful reading  of the paper and a detailed list of comments which improved the exposition. 

%%%%%%%%%%%%%%%%%%%%%%%%%%%%%%%%%%%%%%%%%%%%%%%%%%%%%%%%%%%%%%%%%%%%%%%%%%%%%%%%%%%%%%%%%%%

\verb??\\
{\bf Driss Essouabri}\\
Univ. Lyon,
UJM-Saint-Etienne,\\
CNRS, Institut Camille Jordan UMR 5208,\\
Facult\'e des Sciences et Techniques,\\       
23 rue du Docteur Paul Michelon,\\
F-42023, Saint-Etienne, France\\
{\it E-mail address}: driss.essouabri@univ-st-etienne.fr

\medskip

\verb??\\
{\bf Kohji Matsumoto}\\
Graduate School of Mathematics,\\
Nagoya University,\\
Furo-cho, Chikusa-ku,\\
Nagoya 464-8602, Japan\\
{\it E-mail address}: kohjimat@math.nagoya-u.ac.jp

 \end{document}